\documentclass[12pt]{amsart}
\usepackage{amsmath}
\usepackage{amsfonts}
\usepackage{amsthm}
\usepackage{amssymb}
\usepackage{amscd}
\usepackage{epic}
\usepackage{eepic}

\def\Bbb{\mathbf}
\def\mathrm{}
\newcommand{\labell}[1] {\label{#1}}
\newcommand    {\comment}[1] {}

\newcommand{\mutee}[2]{#2}
\newtheorem{thm}{Theorem}[section]
\newtheorem{prop}[thm]{Proposition}
\newtheorem{lemma}[thm]{Lemma}

\newenvironment{example}{\refstepcounter{thm}\smallskip
\noindent {\bf Example~\thethm}\ }{\hskip\hsize plus0pt minus\hsize
\hbox{$\Box$}\smallskip}

\newenvironment{remark}{\refstepcounter{thm}\smallskip
\noindent {\bf Remark~\thethm}\ }{\hskip\hsize plus0pt minus\hsize
\hbox{$\Box$}\smallskip}

\newcommand{\be}{\begin{equation}}
\newcommand{\ee}{\end{equation}}

\newcommand{\lb}[1]{\label{#1}}
\newcommand{\lbl}[1]{\labell{#1}}

\newcommand{\Ref}[1]{{\rm(\ref{#1})}}  %% Roman declaration added.

\newcommand{\al}{\alpha}
\newcommand{\bet}{\beta}

\newcommand{\de}{\delta}
\newcommand{\De}{\Delta}
\newcommand{\na}{\nabla}
\def \e {{e}}
\newcommand{\ep}{\epsilon}
\newcommand {\cal}[1]  {{\mathcal{#1}}}
\newcommand{\gam}{\gamma}

\newcommand{\om}{\omega}
\newcommand{\Om}{\Omega}
\newcommand{\lam}{\lambda}
\newcommand{\nab}{\nabla}
\newcommand{\up}{\upsilon}

\newcommand{\vphi}{\varphi}

\newcommand{\F}{F}

\newcommand{\wht}{\widehat}
\newcommand{\sk}{SKR }

\newcommand{\kb}{\bar}
\newcommand{\dbar}{\bar{\partial}}
\newcommand{\del}{\partial}

\newcommand{\ra}{\rightarrow}

\def\t{\tau}

\newcommand{\we}{\wedge}

\newcommand{\Lap}{\De}

\title{Scalar curvature and holomorphy potentials}
\author{Gideon Maschler}
\address{Department of Mathematics and Computer Science, Clark University,
Worcester, Massachusetts 01610, U.S.A.}
\email{gmaschler@clarku.edu}
\begin{document}

\begin{abstract}
A holomorphy potential is a complex valued function whose complex gradient, with respect to
some K\"ahler metric, is a holomorphic vector field. Given $k$ holomorphic vector fields on 
a compact complex manifold, form, for a given K\"ahler metric, a product of the following 
type: a function of the scalar curvature multiplied by functions of the holomorphy potentials 
of each of the vector fields. It is shown that the stipulation that such a product be itself 
a holomorphy potential for yet another vector field singles out critical metrics for a particular 
functional. This may be regarded as a generalization of the extremal metric variation of Calabi, 
where $k=0$ and the functional is the square of the $L^2$-norm of the scalar curvature.
The existence question for such metrics is examined in a number of special cases.
Examples are constructed in the case of certain multifactored product manifolds.
For the \sk metrics investigated by Derdzinski and Maschler and residing in
the complex projective space, it is shown that only one type of nontrivial criticality 
holds in dimension three and above.

%An Euler-Lagrange equation is obtained for a family of functionals whose integrands
%involve scalar curvature and appropriately normalized holomorphy potentials, varying over
%K\"ahler metrics associated with a fixed K\"ahler class. This is a simple generalization of
%a variational problem considered by Calabi, which led to his definition of extremal
%K\"ahler metrics. Special cases that are examined include holomorphically equivariant objects
%such as the Futaki invariant, and special K\"ahler-Ricci potentials. A remaining open problem
%is indicated, along with associated iteration problems.
\end{abstract}

\maketitle

\section{Introduction}
\setcounter{equation}{0}

Extremal K\"ahler metrics were defined and studied by Calabi \cite{ca}. In his framework
they are obtained from the variation of a curvature functional, over K\"ahler metrics
(with K\"ahler forms) in a fixed deRham cohomology class of a compact complex manifold.
The resulting Euler-Lagrange equation is equivalent to the requirement that the scalar
curvature of the critical metric is a {\em holomorphy potential}. In other words, the
complex gradient of the scalar curvature of the critical metric is a holomorphic vector field.
%In fact, as the scalar curvature is not
%complex valued but real, the vector field is also Killing, i.e. an infinitesimal generator of isometries.

In this paper we will be concerned with a more general requirement
on the scalar curvature. To put it in context, consider first the stipulation that
\be\lbl{I} \text{the scalar curvature functionally depends
on a collection of holomorphy potentials.}\end{equation}
This condition is satisfied by a subset of the class of
K\"ahler manifolds with rigid torus action, which includes metrics with
Hamiltonian two forms, toric metrics, as well as other classes \cite{ham1}. They are
defined by requiring the values of the metric on vector fields tangent to
a given isometric torus action be functions of the Killing potentials of these vector
fields. Such manifolds are locally classified to be principal torus bundles
over a K\"ahler base, and condition \Ref{I} is satisfied exactly when the
K\"ahler base metric has constant scalar curvature.

The condition we investigate in this article is closely related, but also significantly
more stringent. Namely, we require that
\be\lbl{II} \text{$\vphi_{k+1}:=p(s) h_1(\vphi_1)\cdot\ldots\cdot h_k(\vphi_k)$
is a holomorphy potential},\end{equation}
where $s$ is the scalar curvature, each $\vphi_i$, $i=1\ldots k$ is a holomorphy
potential while $p:\mathbb{R}\ra\mathbb{R}$, $h_i:\mathbb{C}\ra\mathbb{C}$, $i=1\ldots k$
are arbitrary smooth functions.

One reason to concentrate on condition \Ref{II} is that it appears as the Euler-Lagrange
requirement for a functional variation which naturally generalizes the extremal metric
case. In fact, the functional has the form
$$S=S_{f,h_1\ldots h_k}:=\int_{M^m}f(s)h_1(\vphi_1)\cdot\ldots\cdot h_k(\vphi_k)\,\om^{\we m},$$
and its variation gives \Ref{II} for $p=f'$. The special case of extremal metrics
corresponds to the choice $k=0$ and $f(x)=x^2$, i.e. when the functional is the square of the $L^2$-norm
of the scalar curvature. Note that the variation takes place in
a given K\"ahler class while fixing a background metric and $k$ holomorphic vector fields 
(see \S\ref{func}). The functions $f$ and $h_i$, $i=1\ldots k$ are treated as parameters
that must be specified in order to determine the functional.

Returning to condition \Ref{II}, one major distinction is to be made. The vector field
produced by the variation with holomorphy potential $\vphi_{k+1}$ may be linearly dependent on
the $k$ vector fields fixed in defining the variation.
Alternatively, it may be ``genuinely new", i.e. this vector field, together with the fixed ones, 
form a linearly independent set. We explore both cases, focusing mainly on the latter, which is more involved.

In \S\ref{special} we examine for $k=1$ the case of linear dependence, where the constructed vector field
is proportional to the initial fixed one. This case is roughly equivalent to the requirement that the critical metric
satisfying $s=H(\vphi)$, with $\vphi=\vphi_1$ and some smooth $H:\mathbb{C}\ra\mathbb{R}$. A
well understood class of functions fulfilling the latter condition is provided by the class of
\sk metrics \cite{dr-ma1} (see \S\ref{exist} for their definition). These metrics where initially defined
due to their role in the classification of K\"ahler metrics conformal to an Einstein metric (an example
for the latter is the metric constructed by the physicist D. Page \cite{p}).

In \S\ref{exist} we construct, for any $k$, a critical metric with a vector field which forms a linearly independent
set with the fixed ones. This metric is a product of certain \sk metrics. In \S\ref{non-ex} we prove Theorem \ref{thm}, 
whose first part consists of a characterization of \sk metrics among a family of metrics given on a complex vector space, and
an explicit determination of the associated holomorphy potentials of linear vector fields. In the proof of the 
second part of this theorem we examine whether such metrics, aside from their role as critical metrics 
in the linearly dependent case, may also serve as critical metrics (for $k=1$), in the case where the 
constructed vector field and the fixed one are linearly independent. We find a dimensional obstruction 
for this to hold. For some of these \sk metrics the result extends from the vector space to the  
complex projective space.

\section{A variational characterization}\lbl{vary}
\setcounter{equation}{0}
\subsection{Vector fields - real and complex viewpoints}
The material presented in this subsection is roughly equivalent to that in
\cite[Theorem 4.4 and corollaries]{kob}. However, our discussion proceeds mostly
via real, as opposed to complex, geometric terms. Moreover, our aim is to relate the real
and complex viewpoints concretely, and in this we go one small step further than the above reference.

Let $(M,g)$ be a compact K\"ahler manifold with associated almost complex
structure $J:TM\ra TM$. Thus $J^2=-1$, and $J$ is skew-adjoint
and parallel with respect to the Riemannian metric $g$.
The K\"ahler form of $g$ is given by $\om(\cdot,\cdot)=
g(J\cdot,\cdot)$, and the $\mathbb{C}$-valued sesquilinear form whose
real part is $\om$ is $\om^{\mathbb{C}}=(\om-i\om(J\cdot, \cdot))/2=(\om+ig)/2$.
Finally, the induced action of $J$ on $1$-forms is $J^*\xi(\cdot)=\xi(J\cdot)$.

We consider further the action of $J$ on $1$-forms. Let $\Om^k$ denote the space of
$k$-forms on $M$, $k=1\ldots\dim M$. We write $d\Om^k$, $d^*\Om^k$ for the
spaces of exact $(k+1)$-forms and coexact $(k-1)$-forms, respectively, and ${\cal{H}}$
for the space of harmonic $1$-forms. We have
\be\lbl{L2}\text{$J^*(d\Om^0)$ is $L^2$-orthogonal to $d\Om^0$ and to
${\cal{H}}$.}\end{equation}
In fact, by the Hodge decomposition the space $d\Om^0\oplus {\cal{H}}$ is exactly the
space of closed $1$-forms (since for a closed $1$-form, the term of the Hodge 
decomposition in $d^*\Om^2$ is both closed and coclosed, i.e. harmonic and therefore zero). 
Now for a real valued function $\al$ and a closed $1$-form $\nu$, we have
$\int \langle J^*d\al, \nu\rangle  = \int J_p^q \al_{,q} \nu^p = -\int \al\, \omega_{pq} \nu^{p,q} = 0$,
as  $2\omega_{pq} \nu^{p,q} = \omega_{pq} (d\nu)^{qp}$.
%orthogonality to $d\Om^0$ follows from integration by parts, as for
%real valued functions $\al$ and $\bet$ we have
%$\int_M \langle J^*d\bet,  d\al\rangle= -\int_M \al\,\mathrm{Tr}_g(J\circ\nab d\bet)$,
%and this vanishes as $J$ is skew-adjoint while the Hessian of $\bet$ is self-adjoint.
%Next, let $\De_d=-dd^*-d^*d$ be the Laplacian, with $d^*$ denoting the divergence.
%By Bochner's formula, with $\mathrm{r}$ denoting the Ricci tensor,
%\begin{multline*}
%[\De (J^*\xi)]_k=(J^*\xi)_{k,s}^{\ \ \ s}-\mathrm{r}_k^s(J^*\xi)_s
%=J_k^{\ p}\xi_{p,s}^{\ \ \ s}-\mathrm{r}_k^{\ s}J_s^{\ p}\xi_p
%=J_k^p\left(\left(\De\xi\right)_p+\mathrm{r}_p^{\ s}\xi_s\right)
%-\mathrm{r}_k^{\ s}J_s^{\ p}\xi_p\\
%=[J^*(\De\xi)]_k+\left(J_k^{\ l}\mathrm{r}_l^{\ p}-\mathrm{r}_k^{\ s}J_s^{\ p}\right)\xi_p
%=[J^*(\De\xi)]_k,
%\end{multline*}
%the last equality following because the Ricci operator commutes with $J$ on a K\"ahler manifold.
%In particular, of course, $J^*{\cal{H}}\subseteq {\cal{H}}$.
%It follows that for a harmonic $1$-form $\hbar$, we have $\int_M\langle J^*d\bet,\hbar\rangle=
%\int_M\langle d\bet(J\cdot ), \hbar(\cdot)\rangle =-\int_M \langle d\bet,J^*\hbar\rangle$, which
%vanishes as one sees by recalling that $d\Om^0$ is $L^2$ orthogonal to the harmonic forms,
%by the Hodge decomposition.
This completes the verification of \Ref{L2}, which implies, again
by the Hodge decomposition, \be\lbl{sbst} J^*(d\Om^0)\subseteq d^*\Om^2.\end{equation}

%Next, for a closed $1$-form $\eta$, write the Hodge decomposition
%$\eta=d\bet+\hbar$ with $\bet$ a real valued function and $\hbar$ a harmonic $1$-form.
%We will use this below in a setting in which $\eta=\xi-J^*d\al$ for a $1$-form
%$\xi$ and a function $\al$, where we see using \Ref{sbst}, that
%\be\lbl{xi}\text{$d(\xi-J^*d\al)=0 \Longleftrightarrow$ $\xi$ has Hodge decomposition
%$\xi=d\bet+J^*d\al+\hbar$}\end{equation} for $\bet$ and $\hbar$ as above. Equivalently,
%the condition on $\xi$ is that $d\xi$ is the $(1,1)$-form $(i/2)\del\dbar\al$.

Suppose $X$ is a complex vector field of type $(1,0)$ , so that $X=(u-iJu)/2$ for a
vector field $u$, and let $\vphi=\al+i\bet$ be a complex valued function. Then
\be\lbl{C-R} \imath_X\om^{\mathbb{C}}=\dbar\vphi\quad \Longleftrightarrow
\quad \imath_ug=d\bet+J^*d\al,
\end{equation}
where $\imath$ denotes interior multiplication.
%with $J^*\xi(\cdot)=\xi(J\cdot)$ for any one form $\xi$.
This can be seen by substituting for $X$, as well as for
$\dbar\vphi=(d\vphi+iJ^*d\vphi)/2$ with $d\vphi=d\al+id\bet$,
and separating real and imaginary parts. 

Assume $X$ is holomorphic, that is, has holomorphic components in a some complex
coordinate system around each point. Then by \cite[Corollary 4.5]{kob}, it 
satisfies the left hand side of \Ref{C-R} for some function $\vphi$ if it has a nonempty zero set.
Hence the right hand side also holds, and we presently examine $\al$ and $\bet$ in
a special case, namely that of a Killing vector field. The function $\vphi$ is  called a
{\em holomorphy potential} for $X$. Note that for a given holomorphic vector field $X$, the
holomorphy potential varies with the metric, and this will be examined in the next subsection.

Suppose $\t$ is a {\em Killing potential}, i.e a smooth real valued function for which $J\nab\t$
is a Killing vector field. We set $v=\nab\t$ and $u=Jv$. Thus $u$ is an infinitesimal generator
for isometries, that is,  ${\cal{L}}_u g=0$, with ${\cal{L}}$ denoting the Lie derivative.
By \cite[Lemma 5.2]{dr-ma1}, $v$ is holomorphic, that is, ${\cal{L}}_v J = 0$, which is
equivalent to $X=(u-iJu)/2=i(v-iJv)/2$ being holomorphic in the sense of the previous paragraph. 
Also $v$, $u$ (and hence $X$) each have a nonempty zero set consisting of the critical points of 
$\t$. Following the discussion in the previous paragraph, $X$, and hence $u$, satisfy the conditions 
in \Ref{C-R} for some function $\vphi=\al+i\bet$. We wish to describe how $\t$ is related to $\vphi$.

As is well known, $u$ is Killing if and only if $\nab u$ is skew-adjoint. Thus $u$ is
divergence-free, and taking the divergence on the right hand side equation in \Ref{C-R},
gives $$0=\Lap_d\bet+d^* J^*(d\al),$$ with $\Lap_d$ denoting the Laplacian. However,
by \Ref{sbst} $d^*J^*(d\al)=0$, implying that $$\text{$\bet=$ constant}.$$
%Assume now $\xi=\imath_v\omega$ for some vector field $v$. Then $d\xi (\cdot,\cdot)=
%g((\nab u-[\nab u]^*)\cdot,\cdot )$ for $u=Jv$ and $*$ denoting the adjoint (see the proof of
%\cite[Lemma 5.5]{dr-ma1}).
Now set $\xi=\imath_ug$. 
%The right hand side of \Ref{xi} holds (with $\bet$ constant
%and $\hbar=0$, see also \Ref{C-R}). Thus the left hand side also holds, as of course can be
By taking $d$ of the right hand side of \Ref{C-R} we see that $d\xi=dJ^*d\al$. Combining this with
$d\xi(\cdot, \cdot)=2\nab d\t (J\cdot,\cdot)$ (see \cite[Lemma 5.5]{dr-ma1})
gives \be\lbl{bet} \text{$2\nab d\t (J\cdot,\cdot)=d[J^*(d\al )](\cdot,\cdot)=
-2\cdot$\,skew$[\nab d\al (J\cdot,\cdot)]$},\end{equation} where ``skew" denotes the skewsymmetric component.
Note that $\nab d\t$ is Hermitian by \cite[Lemma 5.2 (iii)]{dr-ma1} so the left hand side of
\Ref{bet} is already skewsymmetric. In particular, to represent $X$ via a holomorphy potential $\vphi$
one can always take $$\vphi=\al=-\t.$$ In other words, we can choose the holomorphy potential of
the holomorphic vector field corresponding to a Killing vector field with a nonempty zero set to be,
up to sign, the Killing potential.

Finally, note that a converse also holds. If $X=(u-iJu)/2$ is holomorphic with a {\em real}
holomorphy potential $-\t$, then $u$ is Killing with Killing potential $\t$. In fact, for $v=-Ju$
we have that $-iX=(v-iJv)/2$ is also holomorphic, and this is equivalent to $v$ being holomorphic (in the sense 
that ${\cal{L}}_vJ=0$). Also, $u=Jv$ and $\imath_ug=J^*d\al=-J^*d\t$ (which is the right hand side of \Ref{C-R}
in this case) imply, by composing the latter relation with $J$, that $v=\nab\t$. Hence $v$ is a 
holomorphic gradient, which, by \cite[Lemma 5.2]{dr-ma1} implies that $\t$ is a Killing potential 
for the Killing field $u$.

Due to the above, in the cases below where both notations are used, $\vphi$ and $\t$ will be nearly interchangeable.

\subsection{Functionals}\lbl{func}

In this section we review a lemma of Calabi and then construct a rather direct generalization
to his result on the variation that leads to, and in fact defines, extremal K\"ahler metrics.
We assume throughout that $(M,g)$ is a compact K\"ahler manifold of complex dimension $m$. In the
next lemma we employ complex coordinates $z^a$, with commas denoting complex covariant
differentiation.

\begin{lemma}[Calabi] \lb{l}
If a complex valued function $\vphi$ satisfies the equation
$\vphi^{,ab}_{\ \ ab}=0$, then $\vphi_{,\kb{a}\kb{b}}=0$. In other words,
$\vphi$ is a holomorphy potential (as $\vphi^{,a}_{\ \kb{b}}=0$).
\end{lemma}

\begin{proof}
Write $$0\leq g^{a\kb{c}}g^{b\kb{d}}\vphi_{,\kb{c}\kb{d}}\kb{\vphi}_{,ab}=
(\vphi^{,ab}\kb{\vphi}_{,a}-\vphi^{,ab}_{\ \ a}\kb{\vphi})_{,b}+
\vphi^{,ab}_{\ \ ab}\kb{\vphi}:=
v^{b}_{,b}+\vphi^{,ab}_{\ \ ab}\kb{\vphi},$$
which follows since
$\vphi^{,ab}_{\ \ b}=\vphi^{,ba}_{\ \ b}$. %followed by index renaming
Under the condition of the lemma, the first term on the right hand side is the
divergence of the $(0,1)$ part of the form corresponding to $v$:
$$v^b_{,b}=-*d*(v_{\kb{b}}dz^{\kb{b}}),$$ which vanishes upon integration,
so that the integral of the left hand side also vanishes. That integral
is the square of the $L^2$ norm of  the tensor $\vphi_{,\kb{a}\kb{b}}$, hence the
result follows.
%However, this lemma is
%known to us in another sense, because $\vphi^{,ab}_{\ \ ab}=0$ translates
%via Ricci-Weitzenb\"ock into (the complex version of) Equation (2.10) of
%the main preprint.
\end{proof}
\begin{prop} \lb{S}
On a compact complex manifold $M$ with $\dim_{\mathbb{C}}M=m$, fix objects $\left(\wht{g},\, \{X_i\}_{i=1}^k\right)$,
consisting of a (background) K\"ahler metric $\wht{g}$ with K\"ahler form $\wht{\om}$, and for $i=1,\ldots k$,
a holomorphic vector field (with zeros) $X_i$ admitting a $\wht{g}$-holomorphy potential $\wht{\vphi}_i$.
Suppose that $f:{\mathbb{R}}\ra{\mathbb{R}}$ and $h_i:\mathbb{C}\ra\mathbb{C}$
are smooth functions ($i=1\ldots k$). Consider the functional
$$S=S_{f,h_1,\ldots h_k}:=\int_Mf(s)\prod_{i=1}^k h_i(\vphi_i)\,\om^{\we m},$$ defined over
the space of K\"ahler metrics (with K\"ahler forms) in the fixed K\"ahler class $[\wht{\om}]$. Here,
a given such metric $g$ has K\"ahler form $\om$ and scalar curvature $s$. Also, for i=1\ldots k, each $X_i$ has
$g$-holomorphy potential $\vphi_i$ normalized via the condition
$\int_M\vphi_i\,\om^{\we m}=\int_M[\wht{\vphi}_i+X_iu]\,\om^{\we m}$,
where $u$ is the smooth purely imaginary function given up to a constant by $\om=\wht{\om}+\del\dbar u$.
A metric $g$ satisfies the Euler-Lagrange equation for $S$ if and only if the function
\be\lbl{phi2}\vphi_{k+1}=f'(s)\prod_{i=1}^k h_i(\vphi_i)\end{equation} is
also a $g$-holomorphy potential (and, of course, is a $g$-Killing potential if it is real).
\end{prop}
Here the prime denotes differentiation while $f(s)$ stands for the composition  
$f\circ s$, and similarly for $h_i(\vphi_i)$. Unless all $h_i$, $i=1,\ldots k$ are chosen to
be real valued, Theorem~\ref{S} describes a variation of a complex valued functional.
However, our main interest will be the case where the critical equation \Ref{phi2} gives rise not
just to a holomorphy potential, but to a Killing potential with respect to the critical metric. 
Finally, the meaning of the normalization condition for $\vphi_i$ will be clarified in \S\ref{trivial}.
%Finally, the domain of $f$ will
%turn out to be $R^+$ in one example, requiring an additional assumption on the
%image of $s$.

%\begin{remark}
%Note also that a fixed functional $S$ allows one to consider an iterative process, where
%the pair $(\wht{g},X)$ is replaced by a similar pair consisting of a critical metric
%and its associated "critical" vector field, which now becomes a new input pair for the functional.
%Such a process is well-defined only if there is a procedure for picking one critical metric
%among many.
%\end{remark}
%(along with a fixed based metric $\wht{g}$ allows one to
%consider the following iterative procedure: $X$ is
%regarded as an initial vector field denoted $X_0$ (with some preassigned
%$\wht{g}$-holomorphy potential $\wht{\varphi}$), with a critical metric $g_0$
%distinguishing a second vector field $Y$, denoted $X_1$. Now consider $X_1$
%as an initial vector field, with preassigned potential the one obtained in the previous
%variation stage. A critical metric $g_1$ now gives rise to a new distinguished vector field
%$X_2$. Continuing we get a sequence of critical metrics $g_0, g_1, \ldots$ and vector fields
%$X_0, X_1, \ldots$. One can then investigate conditions on $\{g_i\}$ under which the sequence
%$\{X_i\}$ converges in the (finite dimensional) Lie algebra of vector fields with zeros.
\begin{proof}
We first recall Calabi's basic variation scheme. As the variation takes place in a
fixed K\"ahler class, the K\"ahler forms in a one parameter family $g_t$, may be considered
to vary according to $\om_t=\wht{\om}+\del\dbar u_t$, for purely imaginary functions $u_t$
on the product of $M$ with an open interval.
The $t$-derivative at zero of $u_t$ will be denoted $\up$, and similar derivatives involving
other geometric quantities will be denoted using the symbol $\de$. One has
$\de g_{a\kb{b}}=\up_{,a\kb{b}}$. Next, as variations of determinants involve
traces,  $\de(\det g_{a\kb{b}})=\De \up\det g_{a\kb{b}}$, or
$\de(\om^{\we m})=\De \up (\om^{\we m})$. Here $\De$ is the
$\dbar$-Laplacian, $\De \up=-g^{a\kb{b}}\na_a\na_{\kb{b}}\up$. Consequently the
Ricci tensors $(r_t)_{a\kb{b}}=-\left(\log\det (g_t)\right)_{,a\kb{b}}$
admit the variation $\de r_{a\kb{b}}=(\De \up)_{,a\kb{b}}$.
%=-p_{,\kb{c}}^{\ \,\kb{c}}$.
%which is half the
%$d$-Laplacian, and equal to the $\del$-Laplacian, is given on a function $p$ by

Next, the variation of the scalar curvatures $s_t=\langle r_t,g_t\rangle_t$ for
the one parameter family $g_t$ depends on $r_t$, $g_t$ and the inner products $\langle\cdot,\cdot\rangle_t$.
%(which due to the use of
%complex coordinates is half the scalar curvature in Riemannian geometry).
For this notice that from
$g\cdot g^{-1}=Id$ it follows that $\de (g^{a\kb{c}})g_{b\kb{c}}=
-g^{a\kb{c}}\de(g_{b\kb{c}})$, and so, with $R$ standing for the full curvature tensor,
\begin{eqnarray*}
\de s&=&g^{a\kb{b}}\de r_{a\kb{b}}+\de (g^{c\kb{b}}g^{a\kb{d}}g_{c\kb{d}})
r_{a\kb{b}}=
-\De^2 \up+(1-2)g^{c\kb{b}}g^{a\kb{d}}\de (g_{c\kb{d}}) r_{a\kb{b}}\\
&=&-\De^2 \up- \up^{,a\kb{b}} r_{a\kb{b}}=
-\up_{,a\ \,b}^{\ \,a\ \,b}- \up_{,a}^{\ \,b}r_{b}^{\ \,a}\\
&=&(-\up_{,ab}^{\ \ \,a}+\up_{,d}R^{d\ \ \,a}_{\ \,ab})^{,b}- \up_{,a}^{\ \,b}
r_{b}^{\ \,a}=-\up_{,ab}^{\ \ \,ab}+\up_{,d}^{\ \,b}r^{d}_{\ \,b}+\up_{,d}\,
r^{d\ \,,b}_{\ \,b}- \up_{,a}^{\ \,b}r_{b}^{\ \,a}\\
&=&-\up_{,ab}^{\ \ \,ab}+\up_{,d}s^{,d}\ ,
\end{eqnarray*}
where the symmetry of both the Hessian of $\up$ and $r$ was used in the
second and third lines, the Ricci-Weitzenb\"ock formula in the third line,
and the contracted Bianchi identity in the last line.
%Here the signs and
%all numerical coefficients are quite important since there is a crucial
%cancellation. Also note that in the contracted Ricci identity there is no
%factor of 1/2, which is important in the next steps. Presumably this is
%again because the scalar curvature is defined differently.

We proceed to the variation of the holomorphy potentials $(\vphi_i)_t$ of $g_t$,
for a fixed holomorphic vector field $X_i$, with its $\wht{g}$-holomorphy potential 
denoted by $\wht{\vphi}_i$ as in the statement of the proposition.
For K\"ahler forms $\om_t=\wht{\om}+\del\dbar u_t$, we have
$\imath_{X_i}\om_t=\imath_{X_i}\om+\imath_{X_i}(\del\dbar u_t)=\dbar\wht{\vphi}_i-
\imath_{X_i}(\dbar\del u_t)=\dbar\wht{\vphi}_i+\dbar(\imath_{X_i}\del u_t)=
\dbar(\wht{\vphi}_i+X_iu_t)$,
where we have used the fact that $\imath_{X_i}$ anti-commutes with $\dbar$
when $X_i$ is holomorphic. Hence $(\vphi_i)_t=\wht{\vphi}_i+X_iu_t+(c_i)_t$
for some $t$-dependent constant $(c_i)_t$. But the normalization of the holomorphy potentials
implies $\int_M[\wht{\vphi}_i+X_iu_t+(c_i)_t]\,\om_t^{\we m}=\int_M(\vphi_i)_t\om_t^{\we m}=
\int_M[\wht{\vphi}_i+ X_iu_t]\,\om_t^{\we m}$, so that $\int_M(c_i)_t\,\om_t^m=0$. Hence,
with this normalization, the constants $(c_i)_t$ all vanish and
$(\vphi_i)_t=\wht{\vphi}_i+X_iu_t$. Therefore the variation, i.e. the derivative at $t=0$, is
$$\de {\vphi}_i=X_i\up=\langle\del\vphi_i,\del \up\rangle=
\langle\del\, \mathrm{Re}\, \vphi_i,\del \up\rangle+
i\langle\del\, \mathrm{Im}\, \vphi_i,\del \up\rangle,$$
%\vphi^{,a}\up_{,a}.$$
and the last two terms are, of course, $\de\mathrm{Re}\, {\vphi}_i$ and
$i\de\mathrm{Im}\, {\vphi}_i$, respectively.

Therefore, denoting by $(h_i)_x$, $(h_i)_y$ the partial derivatives
of $h_i$ with respect to $x= \mathrm{Re}\,z$ and $y=\mathrm{Im}\,z$, we have

\begin{eqnarray*}
&&\de \Big(f(s)\prod_{\ell=1}^k h_\ell(\vphi_\ell)\om^{\we m}\Big)=\bigg(f'(s)\left(-\up_{,ab}^{\ \ \,ab}+
\up_{,d}s^{,d}\right)\,\prod_{\ell=1}^k h_\ell(\vphi_\ell)\\
&+&
f(s)\sum_{j=1}^k\Big(\left[\left(h_j\right)_{\! x}\!\!(\vphi_j) \delta(Re\, \vphi_j) +
i\left(h_j\right)_{\! y}\!\!(\vphi_j) \delta(Im\, \vphi_j)\right]
\prod_{\substack{\ell=1 \\ \ell\ne j}}^k h_\ell(\vphi_\ell)\Big)	
+f(s)\prod_{\ell=1}^k h_\ell(\vphi_\ell)\De \up\bigg)\om^{\we m}\\
&=&\bigg(-f'(s)\prod_{\ell=1}^k h_\ell(\vphi_\ell)\up_{,ab}^{\ \ \,ab}\\
&+&\Big\langle f'(s)\prod_{\ell=1}^k h_\ell(\vphi_\ell)\del s+
f(s)\sum_{j=1}^k\Big(\Big(\prod_{\substack{\ell=1 \\ \ell\ne j}}^k h_\ell(\vphi_\ell)\Big)
\left[\left(h_j\right)_{\! x}\!\!(\vphi_j) \del\,(Re\, \vphi_j) +
i\left(h_j\right)_{\! y}\!\!(\vphi_j) \del\,(Im\, \vphi_j)\right]\Big),\del \up\Big\rangle\\
&&+f(s)\prod_{\ell=1}^k h_\ell(\vphi_\ell)\De \up\bigg)\om^{\we m}
=\Big(-f'(s)\prod_{\ell=1}^k h_\ell(\vphi_\ell)\up_{,ab}^{\ \ \,ab}+
{\textrm{div}}(f(s)\prod_{\ell=1}^k h_\ell(\vphi_\ell)\del \up)
\Big)\om^{\we m},
\end{eqnarray*}
%\begin{eqnarray*}
%\de (f(s)h(\vphi)\om^{\we m})&=&\left(f'(s)\de s\,h(\vphi)+
%f(s)h'(\vphi)\de \vphi+f(s)h(\vphi)\De \up\right)\om^{\we m}\\
%&=&\left(f'(s)(-\up_{,ab}^{\ \ \,ab}+\up_{,d}s^{,d})h(\vphi)+
%f(s)h'(\vphi)\vphi^{,a}\up_{,a}+f(s)h(\vphi)\De \up\right)\om^{\we m}\\
%&=&\left(-f'(s)h(\vphi)\up_{,ab}^{\ \ \,ab}+\langle f'(s)h(\vphi)\del s+
%f(s)h'(\vphi)\del\vphi,\del \up\rangle+f(s)h(\vphi)\De \up\right)\om^{\we m}\\
%&=&\left(-f'(s)h(\vphi)\up_{,ab}^{\ \ \,ab}+{\mathrm{div}}(f(s)h(\vphi)\del \up)
%\right)\om^{\we m},
%\end{eqnarray*}

\noindent
using only the Leibniz rule, with ${\textrm{div}}$ denoting the divergence
operator.
Since the second summand is  indeed a divergence, its integral
vanishes. Therefore,
$$\de S=-\int_{M^m}\,f'(s)\cdot\prod_{i=1}^k h_i(\vphi_i)\cdot \up_{,ab}^{\ \ \,ab}
\,\om^{\we m}=
-\int_{M^m}\,\left[(f'(s))\cdot \prod_{i=1}^k h_i(\vphi_i)\right]^{,ba}_{\ \ \,ba}
\cdot \up\,\om^{\we m},$$

\noindent
as one sees integrating by parts four times. But equating this to zero
is a requirement that must hold for every $\up$,  hence one arrives at
$$\left[f'(s)\cdot\prod_{i=1}^k h_i(\vphi_i)\right]^{,ba}_{\ \ \,ba}=0.$$
The result now follows from Lemma~\ref{l}

\end{proof}

\subsection{Remarks on special cases}\lbl{special}
\setcounter{equation}{0}

The following is a description of some examples of metrics with a holomorphy potential
given by \Ref{phi2}.

\subsubsection{Tautological examples leading to invariants}\lbl{trivial}
In some cases Equation \Ref{phi2} gives the holomorphy
potential of the zero vector field. Obviously this can be achieved by
choosing $f$ or one of the functions $h_i$ to be the zero function.
%or setting the initial vector field $X=0$ with $\wht{\vphi}=0$.
Another more interesting possibility is to have $f$ be a nonzero constant function. In that case,
the integrand of $S$ is a product of functions of the $\vphi_i$, while Equation \Ref{phi2} yields
again zero as a holomorphy potential. As this is a holomorphy potential for {\em any} metric,
any metric is critical, hence the functional $S$ must be constant, i.e. an invariant
associated with the K\"ahler class. This well-known fact can be understood in the context
of holomorphic equivariant cohomology. Namely, let $k=1$ and write $X:=X_1$,
$\wht{\vphi}:=\wht{\vphi}_1$, $\vphi:=\vphi_1$ and $h:=h_1$. Denote by $\dbar$ the Dolbeault
operator, and $\imath_X$, as before, interior multiplication by $X$. Then the operator $\dbar-\imath_X$ acts on
equivariant forms, i.e forms in $A_l:=\oplus_{q-p=l}A^{p,q}$, with $A^{p,q}$
denoting $(p,q)$-forms (cf. \cite{l}). Consider for simplicity the case where $h$ is the
identity function, so that the functional is a multiple of $\int_M\vphi\om^{\we m}$.
Then the otherwise mysterious normalization condition in Proposition \ref{S}
can be employed as follows:
$(m+1)\int_M\vphi\,\om^{\we m}=
(m+1)\int_M[\wht{\vphi}+ Xu]\,\om^{\we m}
=(m+1)\int_M[\wht{\vphi}+
Xu]\,\left(\wht{\om}+\del \dbar u\right)^{\we m}
=\int_M[\wht{\om}+\del \dbar u+\wht{\vphi}+Xu]^{\we (m+1)}
=\int_M[\wht{\om}+\wht{\vphi}+\del \dbar u+Xu]^{\we (m+1)}
=\int_M[\wht{\om}+\wht{\vphi}-(\dbar-\imath_X)(\del u)]^{\we (m+1)}$.
%\begin{eqnarray*}
%%\int_M[\om+\vphi]^{\we (m+1)}&=&
%(m+1)\int_M\vphi\,\om^{\we m}&=&
%(m+1)\int_M[\wht{\vphi}+ Xu]\,\om^{\we m}\\
%&=&(m+1)\int_M[\wht{\vphi}+
%Xu]\,\left(\wht{\om}+\del \dbar u\right)^{\we m}
%=\int_M[\wht{\om}+\del \dbar u+\wht{\vphi}+Xu]^{\we (m+1)}\\
%&=&\int_M[\wht{\om}+\wht{\vphi}+\del \dbar u+Xu]^{\we (m+1)}
%=\int_M[\wht{\om}+\wht{\vphi}-(\dbar-\imath_X)(\del u)]^{\we (m+1)}.
%\end{eqnarray*}
In other words, as one varies the K\"ahler form $\om$ in the K\"ahler class $[\wht{\om}]$ (by a change of $u$),
$\vphi$ is constrained to vary in such a manner which results in the closed equivariant form $\wht{\om}+\wht{\vphi}$
changing by an equivariantly exact form, and hence its integral must be constant. This shows the constancy of the
functional from the equivariant coohomology viewpoint.
%The statement is that with this variation, the integral
%of this closed equivariant form, as well as all the integrals with integrands of the form $h\circ\vphi$
%remain unchanged (note that $h$ is not constrained to be polynomial (or exponential) as it does in equivariant
%cohomology).
A discussion of K\"ahler class invariance of related integrals based on just
such a variation of a closed equivariant form appears in \cite[paragraph after (17)]{c}.

If both $f$ and $h$ are the identity functions ($f(x)=x$, $h(z)=z$),
then the Euler-Lagrange equation of $S$ is obeyed by metrics $g$ for which
$\vphi$ is a holomorphy potential. By the definition of $\vphi$,
this holds tautologically for every metric in the variation.
Hence $S$ is again a K\"ahler class invariant.
% provided
%that $\vphi$ satisfied the normalization discussed in (\ref{one}). Note that for this
%$f$ (and nonconstant $\vphi$), there are no solution metrics if $h$ is any non-affine function.
Combining this case with the previous one yields a proof of the
invariance of the Futaki invariant, which, for a holomorphic
vector field $X$ with a nonempty zero set and holomorphy potential $\vphi$, is given by
${\cal F}_{[\om]}(X)=\int_M (s-s_0)\vphi\,\om^n$, with  $s_0$ denoting the
average scalar curvature.
%However, this proof has the limitation of assuming a particular
%normalization on $\vphi$, which is immediately seen to be unnecessary, as the expression for
%${\cal F}_{[\om]}(X)$ is independent of normalization. Associated closed
Related equivariant expressions for these invariants appear in \cite{t,c}.

\subsubsection{The nontrivial case}\lbl{notriv}
We continue with the above notations for the case $k=1$. For either
$X=0$ or $h$ a nonzero constant function,
the Euler-Lagrange requirement is for $f'(s)$ to be a holomorphy
potential. This is perhaps the most direct generalization of Calabi's result,
which is just the case $f(x)=x^2$.

The main situation studied in this work is the case where neither $f$ nor all the functions $h_i$
are constant, at least one of them is not the identity function, and some $X_i\ne 0$.
As mentioned in the introduction, in such a situation it is natural to make a distinction regarding the
result of the variation in Proposition \ref{S}. Namely, we can examine whether or not, for the vector field
$X_{k+1}$ having holomorphy potential $\vphi_{k+1}$ and the fixed vector fields $X_1,\ldots X_k$, the set 
$\{X_1,\ldots X_{k+1}\}$ is {\em linearly independent}.

Consider first for $k=1$, the case where $X_2$ is a multiple of $X_1:=X$.
%In fact, in the first case
%Suppose $h$ is nonconstant, and for the critical metric,
%$(f'\circ s)\cdot (h\circ\vphi)$ is also a holomorphy potential for the
%fixed vector field $X$, which we assume is nontrivial. Equivalently,
This translates into the equation $f'(s)h(\vphi)-a\vphi =\textrm{constant}$
for a constant $a$, i.e. a functional relation between $s$ and $\varphi$,
involving an arbitrarily chosen function $h$.
%As $h$ is nonconstant, one arrives at a case more involved than (\ref{three}),
%and in particular more involved than the extremal metric
%case $s=\varphi$.
Rather than dealing with this equation directly, we consider K\"ahler metrics for which
$s=H(\vphi)$ for some smooth function $H:\mathbb{C}\ra\mathbb{R}$. These can be
fitted into the variational framework we have been studying. In fact, one can choose a
functional $S_{f,h}$ for which a metric satisfying $s=H(\vphi)$, for a fixed function
$H$ as above, is critical. In fact, choose
$f$ to be the exponential function (or any function whose derivative is nowhere zero),
and then define $h(z):=z/f'(H(z))$.
%One natural choice which follows this
%prescription is to take $f$ to be the exponential function.
This yields the functional $$S=\int_M \exp\left[s-H(\vphi)\right]\varphi\,\om^{\we m}$$
whose critical metrics are required to have the integrand as holomorphy potential. If a metric
satisfies $s=H(\vphi)$, that integrand equals $\vphi$, which is of course a holomorphy potential,
so that such a metric is indeed critical, and the resulting ``new" vector field $X_2$ is just $X$.

Note that in the examples in the next paragraph, the critical metric has $\vphi$ as a Killing potential
for $X$ (and hence chosen to be real valued), and not just a holomorphy potential. Hence the equation $s=H(\vphi)$
need only hold for $H:\mathbb{R}\ra\mathbb{R}$. For the purpose of the variation, one can then extend
the domain of $H$ to $\mathbb{C}$ in an arbitrary (smooth) fashion.
% The case which is of most
%interest occurs when for the critical metric, $\varphi$ is in fact a Killing potential.
%In considering just such a metric and not the variation leading to it, one may restrict attention
%from the entire domain of $H$, and consider it instead to be just the real numbers.

As mentioned in the introduction, at least for a large class of smooth functions $H$
(aside from the identity, and with $\vphi$ nonconstant), metrics for which $s=H(\vphi)$, do exist
%on the real numbers, which is arbitrary except for
%satisfying certain boundary conditions, do exist
on compact K\"ahler manifolds. A class of such metrics that will prove quite useful,
is the class of \sk metrics, i.e. metrics admitting a special K\"ahler-Ricci potential,
in complex dimension at least three (see \cite[Lemma 11.1]{dr-ma1} and \cite{dr-ma2}).
A precise definition will be given in \S\ref{exist}. For these metrics, $H$
cannot be completely arbitrary, as it depends on certain boundary conditions (see
again \S\ref{exist}).
%In these examples, the function $\varphi$ is, in fact, a Killing potential.
%In the  special case when these metrics are almost everywhere conformal to
%Einstein metrics, $H$ is, in general, a fairly complicated rational function.

Finally we come to the case studied in the next section, where $\{X_1,\ldots X_{k+1}\}$ is
a linearly independent set.
%in which the
%critical metric has holomorphy potential $\vphi_2=(f'\circ s)\cdot (h\circ\varphi)$
%for a vector field {\em linearly independent} from the initial vector field $X$.
%One remaining intriguing problem is whether there exists a choice of $f$,
%$h$ and nontrivial $X$, for which the resulting functional has a critical metric
%having the property that $(f'\circ s)\cdot (h\circ\varphi)$ is a holomorphy
%potential for a nontrivial vector field $Y$ that is {\em (generically) linearly independent} from $X$.
In \S\ref{exist} we will give a construction of such metrics on certain product manifolds, using products
of \sk metrics, which works for particular functions $f$, $h_i$.
In \S\ref{non-ex} we prove, roughly, for $k=1$ and with a dimensional assumption,
that one can produce no such examples from a certain class of \sk metrics.
%amounts to the condition that a ratio of two distinct (functions of) holomorphy
%potentials, representing (generically) linearly independent vector fields, is a function of the
%scalar curvature. The author is not aware of any known examples on compact manifolds.
%Another  natural place to look further for such examples seems to be the class of K\"ahler metrics
%admitting a Hamiltonian $2$-form \cite{ham1}.
%and in particular, metrics on $CP^m$
%admitting a special K\"ahler-Ricci potential. The main purpose of this note is to
%describe such examples which do not compactify. Along the way a rigidity of gradient holomorphic
%fields is noted.

%While the latter may be taken up in
%future work, note that it leaves open the question of whether such metrics exist that
%are not of the form described in item (\ref{four}).

%At the time this work was first conceived,
%it was asked whether the extremal metric recently found in \cite{clw},
%on the $2$-point blow-up of $CP^2$, had such an additional characterization.
%If this is a case, it could yield a more explicit description of the metric than
%is currently available. However, there is no evidence for this speculation.

\section{Questions of existence}
\setcounter{equation}{0}

\subsection{Existence}\lbl{exist}
\subsubsection{The general construction}\lbl{constr}
Let $(M_i,g_i)$ be a compact K\"ahler manifold with scalar
curvature $s_i$, $i=1,\ldots k+1$. Suppose further that
each $g_i$ admits a nontrivial Killing vector field $X_i$ with Killing potential $\vphi_i$. Assume now that
\be\lbl{conds}\text{$s_i=H_i(\vphi_i)$ for some smooth functions $H_i:\mathbb{C}\ra\mathbb{R}$,
$i=1,\ldots k$ and  $s_{k+1}=\log\vphi_{k+1}$.}\end{equation}
Then the scalar curvature $s$ of the product metric $g=g_1+\ldots +g_{k+1}$
on the product manifold $M=M_1\times\ldots\times M_{k+1}$ satisfies
$s=s_1+\ldots +s_k+s_{k+1}=\sum_{i=1}^k H_i(\vphi_i)+\log\vphi_{k+1}$, so that
%and $e^s=e^{s_1+s_2}=e^{s_1} e^{s_2}=e^{H(\vphi_1)}\vphi_2$ so
\be\lbl{s-exp}\vphi_{k+1}=e^s\prod_{i=1}^k e^{-H_i(\vphi_i)}.\end{equation}
As $\vphi_{k+1}$, and in fact each
$\vphi_i$, may be regarded as a holomorphy potential on $M$ (after pull-back),
it follows that
$$\text{$g$ is critical for $S\!:=\!\!\int_M
\e^{ s- \sum_{i=1}^k H_i(\vphi_i)}\,\om^{\we m}$, \quad $m=\dim_{\mathbb{C}}M$.}$$
In other words, Equation \Ref{phi2} is satisfied and defines a holomorphy potential 
on $M$ for $f(x)=\exp (x)$ and $h_i(z)=\exp (-H_i(z))$. As all vector fields arising 
from different factors of $M$ and are nontrivial, it is clear that $X_{k+1}$ is not 
linearly dependent on the $X_i$, $i=1,\ldots k$.

The question of existence of such examples is thus reduced to that of constructing
the metrics $g_i, i=1,\ldots k+1$. Constructing $g_i$, $i=1,\ldots k$, each
satisfying one of the first $k$ conditions in \Ref{conds} can be achieved, for
many functions $H_i$, using \sk metrics on compact manifolds (with the domain
of $H_i$ enlarged to $\mathbb{C}$ as in \S\ref{notriv}). We complete the
demonstration of existence by showing that $g_{k+1}$ can also be constructed
using such metrics.

\subsubsection{A distinguished \sk metric}
An \sk metric is defined as a K\"ahler metric $g$
which admits a Killing potential $\t$ such that at each $\t$-noncritical
point, all nonzero vectors orthogonal to the complex span of $\nab\t$
are eigenvectors of both the Ricci tensor and the Hessian of $\t$,
considered as operators. The term \sk comes from ``special K\"ahler-Ricci",
a reference to the fact that for such metrics the Hessian of $\t$
satisfies on an open set an equation involving the Ricci tensor and the K\"ahler metric.
The term ``special" refers to the fact that the coefficients in this equation are
not arbitrary, as they are functions of $\t$.

We use a number of results from the classification of \sk metrics \cite{dr-ma1, dr-ma2}.
First, the global classification yields two families on compact manifolds.
One family, on which we will focus in this subsection, is defined on the
projectivization  of a line bundle over a K\"ahler manifold, which is
necessarily Einstein unless it is of complex dimension one.
The other one is defined on $\mathbb{CP}^m$, and will be considered in \S\ref{non-ex}.
%Both families are given by explicit formulas.

To define metrics in the first family, let $(N,h)$ be a K\"ahler manifold, chosen so that
it is also Einstein if $m-1=\dim_{\mathbb{C}} N>1$.
Let $Q$ be a smooth function $Q(\t)$ on a $\t$-interval $[\t_{min},\t_{max}]$,
which is positive on the interior of this interval, and at the endpoints satisfies
boundary conditions (see \cite[Equation (5.1)]{dr-ma2})
\be\lbl{QQ} Q(\t_{min})=Q(\t_{max})=0, \quad Q'(\t_{min})=-Q'(\t_{max})\ne 0.\end{equation}
Choose a complex line bundle $\pi:L\ra N$ with a hermitian metric $\langle\cdot,\cdot \rangle$,
whose associated metric connection has curvature $\Om=-2a\om^{(h)}$, where $\om^{(h)}$ is
the K\"ahler form of $h$ and $a$ is a constant equal to the value of $Q'/2$ at one of the
endpoints. Let $r(\t ):(\t_{min}, \t_{max})\ra (0,\infty)$ be a $C^\infty$ diffeomorphism
satisfying
$dr/d\t=ar/Q$, and fix a constant $c$ smaller than $\t_{min}$. Define a metric on
$\text{$L\setminus\{$ zero section$\}$}$ by
\be\lbl{g} \text{$g=2(\t (r)-c)\pi^*h$ on ${\cal{H}}$, \quad
$g=Q(\t (r))/(ar)^2\,\mathrm{Re}\langle\, ,\rangle$ on ${\cal{V}}$.}\end{equation}
Here ${\cal{V}}$ is the vertical distribution in $TL$ while ${\cal{H}}$ the
horizontal distribution induced by the connection, and ${\cal{V}}$, ${\cal{H}}$
are declared $g$-orthogonal. Also, $\mathrm{Re}\langle\, ,\rangle$ denotes the real part of the hermitian metric while 
$r:L\ra [0,\infty )$ is its norm function, so that $\t=\t(r)$, the inverse of $r(\t)$, becomes a function on $L$.
The metric $g$ and the function $\t$ extend smoothly to the projective compactification of $L$ \cite[Section 5]{dr-ma2}.

As mentioned in the previous section, for an \sk metric in complex dimension at least three, the
scalar curvature $s$ is a function of the Killing potential $\t$ \cite[Lemma 11.1]{dr-ma1}. We will prove
existence of an \sk metric satisfying the condition \be\lbl{s-log}s=\log(\t),\end{equation} by
using it to construct an appropriate function $Q(\t)$, along with an interval $[\t_{min},\t_{max}]$
on which positivity of $Q$  on its interior and the boundary conditions \Ref{QQ} hold. The required metric will be given
via Formula \Ref{g} using this $Q$ and appropriate constants $a$, $c$.

We determine $Q$ using various relationships available to us from
the theory of \sk metrics. Specifically, given any \sk metric, following \cite[Formulas (7.4)]{dr-ma1}, we denote
by $\mu$ and $\lambda$ the eigenvalues of the Ricci tensor on ${\cal{V}}=\mathrm{span}_{\mathbb{C}}(\nab\t)$
and its orthogonal complement $\cal{H}$, respectively. We also denote by $\phi$ the
${\cal{H}}$-eigenvalue of the Hessian of $\t$. As we take the complex dimension $m$ to be at least three,
the functions $\mu$, $\lam$, $\phi$ are functionally dependent on $\t$ (see again \cite[Lemma 11.1]{dr-ma1}).
Finally, the formula $c=\t-Q/(2\phi)$ defines a constant for any \sk metric with $\phi\ne 0$
(see \cite[Lemma 10.1]{dr-ma1}).
The following relations hold among these quantities:
\begin{eqnarray}\lbl{tools}
&(i)&\ s=2\mu+2(m-1)\lam,\\
&(ii)&\ Q\,d\lam/d\t=2(\mu-\lam)\phi\ \textrm{if}\ m>2,\nonumber\\
&(iii)&\ \mu=-(m+1)\phi'-(\t-c)\phi'',\nonumber\\
&(iv)&\ Q=2(\t-c)\phi.\nonumber
\end{eqnarray}
All these relations appear in \cite{dr-ma1}: see Sections $10$ and $11$ for (i),(ii),(iv)
and Section $20$ for (iii).

We now proceed as follows: inserting Equation \Ref{tools}(iv) in \Ref{tools}(ii) one gets, after
cancellation $(\t-c)d\lam/d\t=\mu-\lam$. Replacing $\mu$ in this equation by
$-(m-1)\lam +(\log\t )/2$ (an expression obtained from combining
\Ref{s-log} and \ref{tools}(i)), we get an ODE for $\lambda$, which we solve.
Inserting this solution in the equation obtained by equating
$-(m-1)\lam +(\log\t )/2$  this time with the expression for $\mu$ in \Ref{tools}(iii),
we obtain an ODE for $\phi$. Its solution yields $Q$ using  \Ref{tools}(iv).
To simplify the expressions involved, we only give $Q$ in the case where $m=3$ and
$c=0$:
$$Q=-2\t\, \left((1/24)\,\log  \left( \t \right) \t-(7/288)\,\t+
A\t^{-2}+(1/3)\,B\t^{-3}-C \right)$$ where $A$, $B$ and $C$ are constants.
Next, given interval endpoints $\t_{min}$, $\t_{max}$, the three equations
in \Ref{QQ} give conditions from which $A$, $B$ and $C$ can be determined.
For example, for $\t_{min}=1$, $\t_{max}=2$ we have
$$A={\frac {20}{33}}\,\log  \left( 2 \right) -{\frac {49}{132}},
\quad B=-{\frac {14}{11}}\,\log  \left( 2 \right) +{\frac {115}{132}},
\quad C={\frac 2{11}}\,\log  \left( 2 \right) -{\frac {37}{352}}.$$
As the function $Q$, with the above choices of $A$, $B$ and $C$, satisfies the
boundary conditions, it remains to check whether it is positive in
$(\t_{min},\t_{max})=(1,2)$. We do so using elementary calculus.
First, in this interval $Q>0$ exactly when $\phi>0$ (see \Ref{tools}(iv)).
Writing $\phi=\al \t+\bet \t^{-3}+\gam-\delta \t \ln(\t) -\ep\t^{-2}$,
where all coefficients are positive, we have
$$\t^5\phi''=12\bet-\delta\t^4-6\ep\t,$$ which is decreasing.
One checks that at $\t=1$ this expression, i.e. $12\beta-\delta-6\ep$, is negative,
to conclude that $\phi''<0$, i.e. $\phi'$ is strictly decreasing on
our interval. Checking that $\phi'(1)>0$ and $\phi'(2)<0$, while of course
$\phi(1)=\phi(2)=0$ shows that $\phi'(\t)$ is positive until a unique
critical point in $(1,2)$ on which $\phi$ is positive. For larger $\t$ the function
$\phi$ decreases, so that for $\t$ between the critical point and the zero
at $\t=2$, we also have $\phi(\t)>0$.

Thus $Q$ is shown to be a proper ingredient  for the construction of an \sk metric
satisfying \Ref{s-log}, using Formula \Ref{g}. The other objects needed
are the constant $c$, which is zero in the example above (and so smaller than $\t_{min}=1$),
and the nonzero constant $a$ chosen to equal the value of $Q'/2$ at one of the endpoints, say at $\t =1$.
All together this produces via \Ref{g} a metric on a $\mathbb{CP}^1$-bundle over a 
complex surface (as $m=3$ in our example) admitting a K\"ahler-Einstein metric. This metric
can be taken to be $g_{k+1}$ in the construction leading to Equation \Ref{s-exp}, with $\vphi_{k+1}=\t$.

\subsection{Nonexistence}\lbl{non-ex}
\subsubsection{Prologue}
Let $(M,g)$ be a compact K\"ahler manifold with scalar curvature $s$,
and assume there is a Killing potential $\vphi_0$ for which $s=H(\vphi_0)$
for a smooth function $H:\mathbb{R}\ra\mathbb{R}$ (for instance,
$g$ could be an \sk metric). Suppose one wants to examine whether $g$
{\em also} satisfies the requirement \Ref{II}, namely that $\vphi_{k+1}
:=p(s) h_1(\vphi_1)\cdot\ldots\cdot h_k(\vphi_k)$ is a holomorphy potential,
for holomorphy potentials $\vphi_i$ and some smooth functions $p$ and $h_i$,
$i=1\ldots k$. Now $\vphi_{k+1}$ can be rewritten in the form
\be\lbl{III} \vphi_{k+1}=h_0(\vphi_0) h_1(\vphi_1)\cdot\ldots\cdot h_k(\vphi_k), \end{equation}
where $h_0=p\circ H$. Applying the $\dbar$ operator to both sides of this equation,
and then taking metric duals gives, for the corresponding vector fields $X_i$,
$i=0,\ldots k+1$,  $$X_{k+1}=\sum_{i=0}^k a_iX_i, \quad\quad a_i=h_i'(\vphi_i)
\prod_{\substack{\ell=0 \\ \ell\ne i}}^k (h_\ell(\vphi_\ell)).$$
Now suppose $X_i$, $i=0,\ldots k$ form a basis of the Lie algebra of holomorphic vector fields
with a nonempty zero set, and additionally are pointwise linearly independent somewhere
(and hence, by analyticity, in an open dense set). Then a
necessary condition for $X_{k+1}$ to also be holomorphic is that the $a_i$,
$i=0,\ldots k$ are all constant. This implication places severe restrictions
on the possible functions $h_i$, $i=0,\ldots k$ for which \Ref{III} can be
solved.

Theorem \ref{thm} below describes another situation where Equation \Ref{III}
(or \Ref{II}) has only a restricted type of solution.

\subsubsection{Main theorem}
We consider metrics on a hermitian vector space $(V,\langle\,,\rangle)$ (excluding $0$) of the form
\be\lbl{gg} g=S(r)\mathrm{Re}\langle\,,\rangle_{\cal{V}}+
T(r)\mathrm{Re}\langle\,,\rangle_{\cal{H}}.\end{equation}
for arbitrary smooth and positive functions $S$ and $T$ of $r=|x|$.
Here ${\cal{V}}$ is the distribution in the tangent bundle $TV$, consisting at $x\ne 0$ of
all multiples of $x$, and ${\cal{H}}$ is its orthogonal distribution.
The symbol $\mathrm{Re}\langle\,,\rangle_{\cal{V}}$ denotes the restriction to ${\cal{V}}$
of the real part of the hermitian metric, and similarly for $\mathrm{Re}\langle\,,\rangle_{\cal{H}}$.
%Also, $S$ and $T$ are positive functions of $r:=\sqrt{<x,x>}$.

Among the metrics on  $V \smallsetminus \{0\}$ given by \Ref{gg}, one can find \sk metrics giving rise 
to the second of the two families on compact manifolds mentioned in \S\ref{exist}.
%Both families are given by explicit formulas.
For such metrics \be\lbl{S-T} S(r)=Q(\t (r))/(ar)^2, \quad T(r)=2|\t (r)-c|/(|a|r^2)\end{equation}
(see \cite[Section 6]{dr-ma2}).
Here  $c$ and $a\ne 0$ are constants, and the formula involves a smooth positive
function $Q(\t)$ on a given $\t$-interval $(\t_{min},\t_{max})$. It is required that this interval does not
contain $c$. Given $Q$, the relationship between the potential $\t$ and $r=|x|$ is given again via the equation $dr/d\t=ar/Q$.
As in \S\ref{exist}, this defines $r(\t):(0,\infty)\ra (\t_{min},\t_{max})$, with inverse $\t(r)$.

In Theorem \ref{thm} below, we write $T_r$ for the derivative of $T(r)$. We also use 
the notation \be\lbl{F}\text{$\F=\F_A:=\langle Ax,x\rangle$, with $A$ a {\em self-adjoint} 
linear operator on $V$.}\end{equation}
\begin{thm}\lbl{thm}
%Let $S(r)$, $T(r)$ be positive $C^1$ functions with domain contained in $\mathbf{R}^+$,
Let $g$ be a Riemannian metric on $V \smallsetminus \{0\}$ given by \Ref{gg}, having scalar curvature $s$, and assume
$A$, $F$ are as in \Ref{F}.
\begin{enumerate}
\item The following are equivalent:
\begin{enumerate}
\item $g$ is K\"ahler.
\item $g$ is an \sk metric.
\item The linear ODE $S=T+(r/2)T_r$ holds.
\item All linear vector fields $x\ra Ax$ (with $A$ self-adjoint) are $g$-gradient
with holomorphy potentials given by
\be\lbl{poten} \vphi=(1/2)F\, T + C, \quad C\ \mathrm{constant}.\end{equation}
\end{enumerate}
\item If $\dim V\geq 3$ while any of the conditions in Part $(1)$ holds, and $g$ satisfies the equation
$$\varphi_2=p(s)h(\varphi_1)$$
for holomorphy potentials $\vphi_1$, $\vphi_2$ of linearly independent (self-adjoint) linear vector fields
and smooth functions $p$, $h$, then $h$ is constant.
\end{enumerate}
\end{thm}
\begin{remark}
Part (2) of this theorem states that under certain assumptions, no further nontrivial condition of the form given in \Ref{II}
for $k=1$ holds, with $\vphi_2$ and $\vphi_1$ holomorphy potentials for linearly independent vector fields.
\end{remark}
\begin{remark}
Under certain conditions some of the \sk metrics in Theorem \ref{thm} extend to $\mathbb{CP}^m$, 
with $V$ viewed as embedded in it in the usual way (see \cite[Section 6]{dr-ma2}). Namely,
aside from positivity on a given $\t$-interval $(\t_{min},\t_{max})$, the function
$Q$ is required to satisfy the same boundary conditions as in \Ref{QQ}:
\be\lbl{Q} Q(\t_{min})=Q(\t_{max})=0, \quad Q'(\t_{min})=-Q'(\t_{max})\ne 0.\end{equation}
Furthermore, $a$ must equal the value of $Q'/2$ at one of the endpoints, 
while the constant $c$ must equal one of the endpoint values $\t_{min}$ or $\t_{max}$.

Note that although the underlying space for Theorem \ref{thm} is a (noncompact) vector space of dimension $m$,
for those \sk metrics that extend to $\mathbb{CP}^m$, Part (2) holds on the (compact) 
complex projective space. This follows since the above linear vector fields also
extend to give all gradient holomorphic vector fields on $\mathbb{CP}^m$. 
(Note that as $\mathbb{CP}^m$ is is simply connected, composition with the standard almost complex
structure provides an isomorphism between gradient holomorphic vector fields and Killing fields.)
\end{remark}
The proof of this theorem will be established in a series of propositions.
\begin{prop}\lbl{a-c}
(a) and (c) in Part (1) of Theorem \ref{thm} are equivalent.
\end{prop}
\begin{proof}
Suppose $g$ is a metric as in \Ref{gg}, with K\"ahler form
$\om(\cdot, \cdot )=g(J\cdot,\cdot )$, where $J$ is the standard complex structure
on $V$. Consider
\begin{multline*}
3d\om(v_0,v_1,v_2)=d_{v_0}(\om(v_1,v_2))-d_{v_1}(\om(v_0,v_2))+d_{v_2}(\om(v_0,v_1))\\
-\om([v_0,v_1],v_2)+\om([v_0,v_2],v_1)-\om([v_1,v_2],v_0).\end{multline*}
Using the antisymmetry of $d\om$ one sees that it is enough to check closedness for $\om$ at a point
using triples chosen from a basis of tangent vectors, such that the triple contains no repetition, and
only one choice of order of the vectors needs to be examined. We use a basis induced (generically) by the mutually
orthogonal vector fields $v=rd/dr$, $u=Jv$ and horizontal vector fields (sections of ${\cal{H}}$) $w$,
$w'$, $w''$ which also commute with $v$ and $u$ (in other words, they are projectable as well as horizontal).
The symmetry considerations just mentioned mean that one need only check closedness on the triples $\{u, v, w\}$,
$\{w,w',w''\}$, $\{u,w,w'\}$ and $\{v,w,w'\}$. We now carry out these calculations.

For the first triple, orthogonality as well as the commutation relations imply that
$d\om(v,u,w)=d_w(\om(v,u))=d_w(g(u,u))=d_w(S(r)\mathrm{Re}\langle u,u\rangle)=0$ as
$d_wr=0$.

For the next triple, we use the standard submersion
$\pi:V\setminus \{0\}\ra \mathbb{P}(V)=\mathbb{CP}^{m-1}$.
This map,  with $V$ equipped with $g$ and  $\mathbb{P}(V)$
with the Fubini-Study metric $h$, is a horizontally homothetic
submersion, that is, $\pi^*h=g/f$ for a positive valued function $f$ on
$\pi:V\setminus \{0\}$ whose gradient is vertical, i.e. is in ${\cal{V}}$ (in fact, $f(x)=|x|^2$).
To compute, we denote the K\"ahler form of $h$ by $\om^h$. Also, note that the
horizontal lift operation takes the bracket of two vector fields $w$, $w'$
on $\mathbb{P}(V)$ to $[w,w']^{\cal{H}}$ on
$V\setminus \{0\}$ (where $w$, $w'$ each denote both a vector field on
$\mathbb{P}(V)$ and its horizontal lift with respect to $\pi$ (see \cite[Remark 14.1]{dr-ma1}).
Therefore, a directional derivative term in the expression for $d\om(w,w',w'')$
equals, for example, $d_{w}(\om(w',w''))=d_{w}(f\om^h(w',w''))=f( d_{w}(\om^h(w',w''))$,
due to the verticality of $\nab f$. On the other hand, a term involving the
Lie brackets equals, for example,
$\om([w,w'],w'')=\om([w,w']^{\cal{H}},w'')=f\om^h([w,w'],w'')$. Thus altogether
$d\om(w,w',w'')$ equals $f$ times $d\om^h(w,w',w'')$, which vanishes as
$\om^h$ is closed.

Next, by orthogonality and commutation relations,
$3d\om(u,w,w')=d_u(\om(w,w')-\om([w,w'],u)$. Regard now $V$ 
as the total space of a line bundle over $\mathbb{P}(V)$. The hermitian metric
$\langle\cdot,\cdot\rangle$ induces a hermitian fiber metric on this line
bundle, whose curvature two form is a multiple of $\om^h$. One deduces
from this the relation
\be\lbl{om-h} [w,w']^{\cal{V}}=-2\om^h(w,w')u.\end{equation}
It follows that $\om([w,w'],u)=
g(J[w,w']^{\cal{V}},u)=- 2g(-\om^h(w,w')u,-v)=-2\om^h(w,w')g(u,v)=0$ as $u$, $v$
are $g$-orthogonal. On the other hand $d_u(\om(w,w'))=d_u(g(Jw,w'))=
d_u(T\mathrm{Re} \langle Jw,w'\rangle )=0$ as $u$ is Killing and $d_u r=0$.

Finally, as in the previous case, commutation relations and orthogonality give
$3d\om(v,w,w')=d_v(\om(w,w'))-\om([w,w'],v)$.
The first term on the right hand side is
\begin{eqnarray*}
d_v(g(Jw,w'))&=&r\frac {d}{dr}(T\mathrm{Re}\langle iw,w'\rangle_{\cal{H}})\\
&=& rT_r\mathrm{Re}\langle iw,w'\rangle+T(\mathrm{Re}\langle i\na_v w,w'\rangle+
\mathrm{Re}\langle iw,\na_v w'\rangle )\\
&=&rT_r\mathrm{Re}\langle iw,w'\rangle+
T(\mathrm{Re}\langle \na_{iw} v,w'\rangle+\mathrm{Re}\langle iw,\na_{w'} v\rangle )\\
&=&(rT_r+2T)\mathrm{Re}\langle iw,w'\rangle
\end{eqnarray*}
where we have used in the last line
the fact that $\nab v$, treated as a bilinear form,
is the Hessian of $r^2/2$, which is of course $\mathrm{Re}\langle\cdot ,\cdot\rangle$.

As for the second term in right hand side, we again use relation \Ref{om-h}
to give
\begin{eqnarray*}
\om([w,w'],v)&=&g(J[w,w']^{\cal{V}},v)=- 2g(-\om^h(w,w')u,u)=2\om^h(w,w')g(u,u)\\
&=&2\om^h(w,w')S\mathrm{Re}\langle u,u\rangle_{\cal{V}}
=2[\mathrm{Re}\langle Jw,w'\rangle /\mathrm{Re}\langle v,v\rangle ]S\mathrm{Re}\langle u,u\rangle=2S\mathrm{Re}\langle iw,w'\rangle,
\end{eqnarray*}
where in the second to last equality we have used the definition of the Fubini-Study K\"ahler form, and in the last line we used
the equality $\mathrm{Re}\langle v,v\rangle=\mathrm{Re}\langle u,u\rangle$.

Thus $d\om (v,w,w')=0$ if and only if $rT_r+2T=2S$, i.e. the equation in  Part (1)(c) of Theorem \ref{thm} holds.
This completes the proof of equivalence.
\end{proof}

Next, we have
\begin{prop}\lbl{b-c}
Conditions (b) and (c) of Part (1) of Theorem \ref{thm} are equivalent.
\end{prop}
\begin{proof}
An \sk metric is in particular K\"ahler, so satisfies (c) by Proposition \ref{a-c}.
For the converse, suppose (c) holds for some smooth and positive functions $S(r)$, $T(r)$.
Fix constants $a>0$, $c$ and define $\t(r):=ar^2T/2+c$ and then $Q:=ar\t'(r)$.
Note that $\t>c$ for $r>0$. We compute $T+(r/2)T_r=2(\t-c)/(ar^2)+(r/2)[-4(\t-c)r^{-3}/a+2r^{-2}\t'/a]
=\t'/(ar)=Q/(ar)^2$, so that $S=Q/(ar)^2$, hence $S$ and $T$ are of the form given in Equation
\Ref{S-T} for the coefficients of an \sk metric. Also, $\t'=arT+ar^2T_r/2=ar(T+(r/2)T_r)=
arS>0$, so that $Q>0$ (if $r>0$) as required for an \sk metric. Finally, by Proposition \ref{a-c},
the metric $g$ given by \Ref{gg} (with \Ref{S-T}) is K\"ahler. Hence $g$ is an \sk metric.
\end{proof}

In the course of proving that (c) and (d) in Part $(1)$ of Theorem \ref{thm} are equivalent,
we consider a more general family of metrics. Here $T_r$ will denote a partial derivative of
a function $T(r,\F)$.
\begin{prop}\lbl{grad}
Let S:=S(r,\F), T:=T(r,\F) be positive $C^1$ functions on a domain in  $\mathbb{R}^2$.
%\be\lbl{PDE} 2\F(S-T)\E_\F=E_r+C', \quad \E=\int T d\F, \end{equation}
Consider a metric defined on $V\smallsetminus\{0\}$ for a
vector space $V$, of the form
$$\lbl{g1} g=S(r,\F)\mathrm{Re}\langle\,,\rangle_{\cal{V}}+
T(r,\F)\mathrm{Re}\langle\,,\rangle_{\cal{H}},$$
with $r=|x|$ and  $\F$ (and $A$) as in \Ref{F} but with $A$ not a multiple of the identity.
The linear PDE \be\lbl{PDE} 2[(S-T)\F]_{\F}=r \,T_r \end{equation} holds
if and only if the linear vector field $x\ra Ax$  is gradient with respect to $g$ and
has holomorphy potential
\be\lbl{potgen}\varphi=(1/2)\int T d\F+C(r),
\quad\text{C(r) an arbitrary $C^1$ function.}\end{equation}
\end{prop}
Note that if the metric dependence on $F$ is nontrivial, the above vector field is the only
gradient linear one.
\begin{proof}
Consider the orthogonal splitting $Ax=(\F/r^2)x+y$, with $y\in {\cal{H}}$. The one-form
$P$ dual to $x\ra Ax$ is given by
\begin{eqnarray*}
P=g(Ax,\cdot)&=&(\F r^{-2})S\cdot\mathrm{Re}\langle x,\cdot\rangle +T\cdot\mathrm{Re}\langle Ax-(F r^{-2})x,\cdot\rangle\\
             &=&\F r^{-2}(S-T)\mathrm{Re}\langle x,\cdot\rangle+T\cdot\mathrm{Re}\langle Ax,\cdot\rangle,\\
\end{eqnarray*}
so that \be\lbl{P} P= r^{-1}F(S-T)\, dr+(T/2)\, d\F,\end{equation}
where we have used the easily verifiable relations %(by computing in a basis)
$$\mathrm{Re}\langle x,\cdot\rangle=r\,dr,\quad \mathrm{Re}\langle Ax,\cdot\rangle=(1/2)d\F,$$
with the latter relation holding precisely since $A$ is self-adjoint.
We compute
$$dP=r^{-1}[F(S-T)]_{\F} \,d\F\we dr+(1/2)\, T_r \,dr\we d\F.$$
As $A$ is not a multiple of the identity, this expression vanishes exactly when \Ref{PDE} holds.
Writing $P=d\varphi=\varphi_r\,dr+\varphi_{\F}\,d\F$ together with \Ref{P} gives
\Ref{potgen} after integrating $\varphi_{\F}$ with respect to $\F$. Note that one can also arrive at the
formula for the holomorphy potential from \Ref{PDE}. Namely, the coefficient of $dr$ in \Ref{P} is, up to
the factor $2/r$ just the $\F$-integral of the left hand side of \Ref{PDE}, so that one can equate
$\varphi_r$ with $2/r$ times the $\F$ integral of the right hand side of \Ref{PDE} and then
integrate with respect to $r$ while switching the order of integration (with respect to $\F$ and $r$).
\end{proof}
The proof of the equivalence of (c) and (d) in Part $(1)$ of Theorem \ref{thm} now follows at once, for
linear vector fields for which $A$ is not a multiple of the identity. In fact,
the linear ODE in (c) follows immediately from \Ref{PDE} and the independence from $\F$
of $S$ and $T$. As $\F$ does not appear in this ODE, the conclusion that $x\ra Ax$ is
gradient holds for {\em any} linear vector field with $A$ as in \Ref{F} (but not a multiple of the identity). 
The expression \Ref{poten} for the potential is a consequence of \Ref{potgen} and has $C(r)=C$ constant, 
as one can see by using the ODE in (c) in the equation relating $\vphi_r$ with the coefficient of $dr$ 
in \Ref{P}, and then integrating with respect to $r$.

To complete the proof that (c) and (d) of Part (1) of the theorem are equivalent, it is 
enough to show that (c) implies formula \Ref{poten} for 
the holomorphy potential, in the case where $A$ {\em is} a multiple of the identity.
In fact, a computation as the one leading to \Ref{P} and the material past it shows 
in this case that $d\vphi=P=\vphi_r\, dr$, so that $\vphi$ is a function of $r$. Now the special
K\"ahler-Ricci potential $\t=\t(r)$ is also a holomorphy potential which is a function
of $r$. But the holomorhpic vector fields associated with these two potentials must be
linearly dependent. This can be seen from equating the term $\dbar r$ in the relation
$\imath_X\om=\dbar\vphi=\dbar (\vphi(r))=\vphi'(r)\dbar r$ with the same term in the analogous relation 
for $\t (r)$ (with $\om$ denoting the K\"ahler form of $g$). This gives, after taking the dual,
that the corresponding $(1,0)$-vector fields are pointwise proportional, so their ratio is real-valued, 
yet must be a holomorphic function, hence is constant. Thus $v=\nab\t$ is a linear vector field
for which the corresponding matrix $A$ is a multiple of the identity. But for $\t$ we have, as 
(b) implies (c) in Part (1), that (with $I$ denoting the identity operator) 
$\t=ar^2T/2+c=\langle aIx,x\rangle T/2+c=(1/2)FT/2+c$ for constants $a$ and $c$ 
(see the proof of Proposition \ref{b-c}). In other words, formula \Ref{poten} holds for $\t$ and 
hence also for any other holomorphy potential for a linear vector field for which $A$ is a 
multiple of the identity.

To continue with Part $(2)$ of Theorem \ref{thm}, we need the following lemma, in
which $I$ again stands for the identity operator on a vector space $V$.
\begin{lemma}\lbl{depend}
Let $A$ and $B$ be two self-adjoint operators on a vector space $V$ with $\dim (V)\geq 3$.
If $y$, $Ay$ and $By$ are linearly dependent for
each $y\in V$, then $A$, $B$, and $I$ are linearly dependent in $Hom(V,V)$.
\end{lemma}
\begin{proof}
Clearly if $A$ or $B$ are a multiple of $I$, the conclusion holds. Hence we assume
$A$ is not a multiple of the identity. Choose eigenvectors $u$, $v$ of $A$ with
distinct eigenvalues $\lam$, $\mu$, respectively. For any scalar $t\ne 0$, the vectors
$u+tv$, $A(u+tv)=\lam u+t\mu v$, $B(u+tv)$ are linearly dependent by assumption, while the
first two of these vectors are linearly independent (their coefficient determinant
$t(\mu-\lam)$ is nonzero). Hence, by the assumption of the lemma, for some scalars $a(t)$, $b(t)$
we have
$Bu+tBv=a(t)(u+tv)+b(t)(\lam u+t\mu v)=(a(t)+b(t)\lam )u+t(a(t)+b(t)\mu )v$.
Differentiating with respect to $t$
gives \be\lbl{v,u} Bv=\al (t)u+\bet(t)v\end{equation}
for appropriate coefficients $\al(t)$, $\bet(t)$.

We now make the following observation:
{\em An $A$-eigenvector $v$ which is not a $B$-eigenvector belongs to
an $A$-eigenspace of codimension at most one.} In fact, otherwise
let $u$, $w$ be two linearly independent $A$-eigenvectors not
in the $A$-eigenspace of $v$. The two pairs $\{ v,u\}$ and $\{ v,w\}$ both satisfy an
equation of the form  \Ref{v,u}. If $v$ is not a $B$-eigenvector then the coefficient
$\al (t)$ and the corresponding one for $w$ are both nonzero. Equating the right hand sides
of the two equations together implies that $\{ u,v,w\}$ are linearly dependent. This is
a contradiction to the fact that these vectors do not all belong to the same $A$-eigenspace, 
and the two of them that may be in the same $A$-eigenspace are linearly independent (note that the coefficient 
of $v$ in the linear dependence equation cannot be zero or else $w$ and $u$ will be proportional).

From the observation above it follows that {\em if $\dim (V)\geq 3$, every $A$-eigenvector is
also a $B$ eigenvector}. For suppose not. Then, by the previous paragraph, there is a codimension
one $A$-eigenspace which is not a $B$-eigenspace (recall that $A$ is not a multiple of the identity,
so the codimension is in fact equal to one). Consider the one dimensional $A$-eigenspace orthogonal
to this codimension one eigenspace (they are orthogonal as $A$ is self-adjoint). It cannot be also a
$B$-eigenspace, for then by orthogonality the codimension one $A$-eigenspace will also be a $B$-eigenspace,
contradicting our assumption. Thus, according to the previous paragraph, this one dimensional 
$A$-eigenspace must be of codimension one. This would make $\dim(V)\leq 2$, a contradiction.

Hence, let $\{v_i, i=1,\ldots \dim (V)\}$ be a basis of eigenvectors
for both $A$ and $B$, with corresponding eigenvalues $\lambda_i^A, \lam_i^B, i=1,\ldots \dim (V)$,
given in their multiplicity. Consider $v=\sum_{i=1}^{\dim (V)} v_i$.
Assume $(aA+bB+cI)v=0$ for some fixed scalars $a$, $b$, $c$ not all zero. This is equivalent to
$\sum_{i=1}^{\dim (V)} (a\lam_i^A+b\lam_i^B+c)v_i=0$, or, by the linear independence of the $v_i$,
to a system of linear equations of the form $$a\lam_i^A+b\lam_i^B+c=0, \quad i=1\ldots \dim (V),$$
which is simply equivalent to the diagonal matrix equation $a\langle\lam_1^A,\ldots \lam_{\dim (V)}^A\rangle+
b\langle\lam_1^B,\ldots \lam_{\dim (V)}^B\rangle+cI=0$, so $aA+bB+cI=0.$

%The statement is obvious in dimension one. Assume $\dim (V)=2$, and let $u$, $v$ be as in the
%first paragraph. If $\al (t)=0$ in \Ref{v,u} then the two equations $(aA+bB+cI)v=0$, $(aA+bB+cI)u=0$
%gives the result just as in the previous paragraph. Otherwise from these two equations we have
%$a\lam v+b\al(t) u+b\beta(t) v +cv=0$ and $a\mu u+bBu+cu=0$. These in turn imply, using the linear
%independence of $u$ and $v$ that $b=0$, $a\lam+c=a\mu+c=0$, and the latter force $a=c=0$.
%Hence the assumption of the proposition cannot hold, and it is true vacuously.
\end{proof}

We now consider Part (2) of Theorem \ref{thm}. We need to determine whether, for a metric $g$
given by \Ref{gg} and satisfying $S=T+(r/2)T_r$, two  (not affinely related) holomorphy potentials
$\varphi_1$, $\varphi_2$ (of linear vector fields) may be related by $$\varphi_2=p(s)h(\varphi_1),$$ where $s$
is the scalar curvature of $g$ and $p$, $h$ are arbitrary smooth functions.
Since metrics as in \Ref{gg} are of cohomogeneity one, it is known that the scalar curvature
is a function of the norm $r$ (this can also be seen as $g$ is an \sk metric, hence $s$ is
a function of $\t$, which is a function of $r$). It is thus sufficient
to replace p(s) by a function $u(r)$, giving
\be\lbl{s-phi} \varphi_2=u(r)h(\varphi_1).\end{equation}
Substituting\Ref{poten}, Equation \Ref{s-phi} takes the form
$$T(r)\langle \overline{A}_2x,x\rangle+C=u(r)h\left(T(r)\langle \overline{A}_1x,x\rangle+C_1\right),$$
where the indices correspond to those of the potentials, $\overline{A}_i=A_i/2$, $i=1,2$ and $C_1$, $C$ are constants.

In order to understand the possible solutions to this equation, we express it a bit differently.
Making the change of variable $y=(T(r))^{1/2}x$, we have $|y|^2=r^2T(r)$ and so
\be\lbl{change}\langle \overline{A}_2y,y\rangle+C=
q(|y|^2)\wht{h}(\langle \overline{A}_1y,y\rangle),\end{equation}
with $q(|y|^2)=u(r)$ and $\wht{h}(w)=h(w+C_1)$.

We now have
\begin{prop} If $\dim (V)\geq 3$ then
any self adjoint solutions $\overline{A}_1$, $\overline{A}_2$ of Equation \Ref{change}, together with
the identity element $I$, are linearly dependent.
\end{prop}
\begin{proof}
Differentiating \Ref{change} shows that the gradients of the three quadratic forms appearing
in \Ref{change} are pointwise linearly dependent, i.e. for all $y\in V$
$$\overline{A}_2y=aIy+b\overline{A}_1y$$ for functions $a=q'(|y|^2)\wht{h}(\overline{F}_1)$,
$b=\wht{h}'(\overline{F}_1)q(|y|^2)$, with $\overline{F}_1=\langle \overline{A}_1y,y\rangle$.
We now use Lemma \ref{depend}.
\end{proof}
Suppose $a_0I+a_1\overline{A}_1+a_2\overline{A}_2=0$ for constants $a_i$, $i=0,1,2$ not all zero. If $a_2=0$
then $A_1$ is a multiple of the identity, so that Equation \Ref{change} has the form
$\langle\overline{A}_2y,y\rangle+C=q(|y|^2)\wht{h}(m|y|^2)$ for some constant $m$. Setting $t=|y|^2$
and $z=\langle \overline{A}_2y,y\rangle$ we have $z+C=q(t)\wht{h}(t)$, which implies, via separation of variables, that
$z$ is constant. But then $\overline{A}_2=0$ so $\vphi_2$ is a holomorphy potential associated with the zero vector
field, which contradicts linear independence of the vector fields. Hence $a_2\ne 0$, and we write
$\overline{A}_2=kI+l\overline{A}_1$ for constants $k$, $l$. We now have from \Ref{change}
$$k|y|^2+l\langle \overline{A}_1y,y\rangle+C=
q(|y|^2)\wht{h}(\langle \overline{A}_1y,y\rangle).$$
Setting $t=\langle y,y\rangle$, $z=\overline{F}_1=\langle \overline{A}_1y,y\rangle$ we have
$kt+lz=q(t)\wht{h}(z)$. taking the partial derivatives with respect to $t$ and $z$
we get $$k=q'(t)\wht{h}(z), \quad l=q(t)\wht{h}'(z).$$ Now $k\ne 0$ as $\overline{A}_2$
and $\overline{A}_1$ represent linearly independent vector fields. If $l\ne 0$,
separation of variables now gives that both $q$ and $\wht{h}$ must be constant, but also that
$q'$ and $\wht{h}'$ are constant. Thus the latter two functions must be identically zero, which
contradicts $k$ and $l$ being nonzero. Hence we must have ($k\ne 0$ and) $l=0$. In this case, as $q$ cannot be zero (since it forces
$q'=0$ and hence $k=0$), we must have $\wht{h}'=0$, i.e. $\wht{h}$ is constant, and so is $h$, as required.
(Note that in this case $q'$ is a nonzero constant, so $q$ is affine in $t$. Finally note that in this
case $\overline{A}_2$ is special, being a multiple of the identity, which means that the corresponding vector 
field is a multiple of $\nab\t$, with $\t$ the special K\"ahler-Ricci potential. The resulting
equation $\text{$\vphi_2=$constant$\cdot p(s)$}$ corresponds to the characterization of the scalar curvature 
of an \sk metric as a function of $\t$.) 
This completes the proof of Theorem \ref{thm}.

\vspace{.3in}
\begin{center}
\textbf{Acknowledgements}
\end{center}
The author thanks Andrzej Derdzinski for many discussions that significantly improved
both the style and the mathematical content. In particular, the idea
behind \S\ref{constr} and Lemma \ref{depend} with its application are due to him.

%One of his comments was the germ that led to the introduction of normalizations,
%and the others improved the overall style.

% ************************************************************************
%       Bibliography:
% ************************************************************************


\begin{thebibliography}{99}

\bibitem{ham1} Apostolov, V., Calderbank, D.M.J., Gauduchon P.: Hamiltonian
2-forms in K\"ahler geometry, I. General theory.  J. Differential Geom.  73
(2006), 359--412
\bibitem{ca} Calabi, E.: Extremal K\"ahler metrics.  Seminar on Differential
Geometry,  pp. 259--290, Ann. of Math. Stud., 102, Princeton Univ. Press,
Princeton, N.J., 1982
\bibitem{dr-ma1} Derdzinski, A., Maschler, G.: Local classification
of conformally-Einstein K\"ahler metrics in higher dimensions. Proc. London
Math. Soc. 87 (2003), 779--819
\bibitem{dr-ma2} Derdzinski, A., Maschler, G.: Special K\"ahler-Ricci
potentials on compact K\"ahler manifolds. J. reine angew. Math. 593 (2006), 73--
116
\bibitem{kob} Kobayashi, S.: Transformation groups in differential geometry.
Reprint of the 1972 edition. Classics in Mathematics. Springer-Verlag, Berlin, 1995
\bibitem{l} Liu, K.: Holomorphic equivariant cohomology. Math. Ann.  303
(1995), 125--148
\bibitem{c} Maschler, G.: Central K\"ahler metrics.  Trans. Amer. Math. Soc.  355
(2003), 2161--2182
\bibitem{p} Page, D.: A compact rotating gravitational instanton. 
Phys. Lett. B 79 (1978), 235--238
\bibitem{t} Tian, G.: K\"ahler-Einstein metrics on algebraic manifolds.
Transcendental methods in algebraic geometry (Cetraro, 1994),  143--185,
Lecture Notes in Math., 1646, Springer, Berlin, 1996
\end{thebibliography}
\end{document}